\documentstyle[12pt]{article} 
\input{epsf}
\def\examples{{G}}
\newcommand{\lra}{\longrightarrow}
\newcommand{\complex}{{\bf C}}
\newcommand{\C}{{\bf C}}
\newcommand{\proj}{{\bf P}}
\newcommand{\integer}{{\bf Z}}
\newcommand{\cC}{{\cal C}}

\newcommand{\cF}{{\cal F}}

\newcommand{\cO}{{\cal O}}
\newcommand{\cT}{{\cal T}}
\newcommand{\jb}{{\overline{J}}}
 \newcommand{\IP}{{\bf P}}

\newcommand{\kf}{{\cal F}}

\def\qed{\hfill $\Diamond$}
\def\Def{D}
\def\loc{\hbox{\sl loc}}
\def\cdia#1#2#3#4{\begin{array}{ccc}
#1 & \lra & #2\\
\downarrow & & \downarrow \\
#3 & \lra & #4\\
\end{array}}

\def\section#1{\par\medskip\par{\large\bf #1}\par\smallskip}

\begin{document}
\centerline{\Large \bf Euler number of the compactified Jacobian}
\vskip 10pt
\centerline{\Large \bf and multiplicity of rational curves.}
\bigskip
\centerline{\bf by B. Fantechi, L. G\"ottsche, D. van Straten}
\smallskip
\bigskip
\begin{abstract}In this paper we show that the Euler number
of the compactified Jacobian of a rational curve $C$ with locally
planar singularities is equal
to the multiplicity of the $\delta$-constant stratum in the base of a 
semi-universal deformation of $C$. In particular, 
the multiplicity assigned by Yau, Zaslow and Beauville to a rational
curve on a K3 surface $S$ coincides with the multiplicity of the 
normalisation map 
in the moduli
space of stable maps to $S$.\end{abstract}
\vskip 15pt

\section{Introduction}
Let $C$ be a reduced and irreducible projective curve with singular set
$\Sigma \subset C$ and let $n: \widetilde{C} \lra C$ be
its normalisation. The generalised Jacobian $JC$ of $C$ is an extension of
$J\widetilde{C}$ by an affine commutative group of dimension 
$$\delta:=\dim H^0(n_*(\cO_{\widetilde{C}})/\cO_C)=\sum_{p \in \Sigma} 
\delta(C,p)$$
so that $\dim JC=\dim J\widetilde{C} +\delta=g(\widetilde{C})+\delta$ 
is equal to
the arithmetic genus $g_a(C)$ of $C$. The non-compact space $JC$ 
is naturally an 
open subset of the {\em compactified Jacobian} $\jb C$ of $C$, 
whose points 
correspond to isomorphism classes of rank one torsion free sheaves 
$\cF$ of degree 
zero (i.e. $\chi(\cF)=1-g_a(C)$) on $C$. The space $\jb C$ is irreducible 
if and only 
if $C$ has planar singularities; then $\jb C$ is in fact a 
compactification of 
$JC$, i.e., $JC$ is dense in $\jb C$. If moreover $C$ is rational 
and unibranch, 
then $\jb C$ is topologically the product of compact spaces 
$M(C,p)$ for every 
$p\in \Sigma$.  The space $M(C,p)$ only depends on the 
analytic singularity 
$(C,p)$; it can be defined as $\jb D$ for any rational curve $D$ 
having $(C,p)$ as 
unique singularity.\\
Let $B=B(C,p)$ be the base of a semi-universal deformation of 
the singularity 
$(C,p)$. Inside $B$ let $B^\delta=B^\delta(C,p)$ be
the locus of points  for which $\delta$ remains constant. This 
means that
$$t \in B^\delta \Leftrightarrow \sum_{p \in C_t} \delta(C_t,p) 
=\delta(C).$$
The codimension of $B^\delta$ is $\delta(C,p)$; its multiplicity 
$m(C,p)$ at 
$[(C,p)]$ is by definition equal to the 
number of intersection points with a generic $\delta$-dimensional 
smooth subspace 
of $B$. The $\delta$-constant stratum can be defined in a similar 
way for a 
semi-universal deformation of a projective curve with only planar 
singularities.
In this paper we show the following theorem.

{\bf Theorem 1.}
{\sl Let $(C,p)$ be a reduced plane curve singularity. Then the 
Euler number of 
$M(C,p)$ is
equal to the multiplicity of the $\delta$-constant stratum:
$$e(M(C,p))=m(C,p).$$
Let $C$ be a projective, reduced rational curve with only planar 
singularities. 
Then $e(\jb C)= m(C)$, the multiplicity of the $\delta$-constant 
stratum 
$B^\delta$ at $0$.}

Note that this gives an independent proof of the following result 
of Beauville: 
Let $C$ be an irreducible and reduced rational curve with planar 
singularities. 
Then $e(\jb C)$
can be written as a product over the singularities of $C$ of a 
number only 
depending on the type of the singularity, and it is the same 
for $C$ and its 
minimal unibranch partial normalisation.

Theorem 1 has an application in the following situation.
Let $X$ be a (smooth) $K3$ surface with a complete 
(hence $g$-dimensional) linear 
system of curves of genus $g$.
Under the assumption that all curves in the system are 
irreducible and reduced,
it was shown in \cite{Y-Z} and \cite{B} that the ``number'' 
$n(g)$ of rational 
curves
occuring in the linear system, is equal to the $g^{\rm th}$ 
coefficient of the 
$24^{\rm th}$ power of
the partition function, i.e:
$$\sum_{g \ge 0}n(g)q^g=\frac{q}{\Delta(q)}$$
where $\Delta(q)=q\prod_{n \ge 1}(1-q^n)^{24}$.
In this counting, a rational curve $C$ in the linear 
system contributes 
$e(\jb C)$ to $n(g)$:
$$n(g)=\sum_{C}  e(\jb C).$$
If $C$ is a rational curve with only nodes as singularities, 
then $e(\jb C) =1$, 
so that $e( \jb C)$
seems to be a reasonable notion of multiplicity. 
Theorem 1 implies that $e(\jb C)$ 
is always positive, and in principle allows an 
explicit computation of it (see 
section \examples).\\

In fact, we prove a more precise statement. For any 
projective scheme $Y$ and 
$d \in H_2(Y,\integer)$ let $M_{0,0}(Y, d)$ be the 
moduli space of
genus zero stable maps $f: {\bf P}^1 \lra Y$ with 
$f_*([{\bf P}^1])=d$. 
Under the above assumptions on the K3 surface $X$ 
and the linear system 
corresponding to
$d$, the space $M_{0,0}(X,d)$ is a zero-dimensional scheme.
If $C \stackrel{i}{\hookrightarrow} X$ is a rational 
curve in $X$ (always assumed 
to be irreducible and reduced), $n:\proj^1 \lra C$ its 
normalisation, then $f=i 
\circ n:\proj^1 \lra X$ is a point of
$M_{0,0}(X,d)$. The moduli space 
$M_{0,0}(X,d)$ contains naturally as a closed subscheme 
$M_{0,0}(C,[C])$, the submoduli space of maps whose scheme 
theoretic image
is $C$; the latter scheme is of course defined for any 
projective reduced curve 
$C$, and it is zero-dimensional if the curve is rational.
More generally, $M_{g,0}(C,[C])$ is zero dimensional, where
$g$ denotes the genus of the normalisation of $C$. The following 
theorem gives another 
interpretation of $e(\jb C)$ in terms of the length of such 
zero-dimensional 
schemes.\\

{\bf Theorem 2.}
{\sl Let $C$ be a reduced, irreducible projective curve with 
only planar  
singularities, and let $g$ be the genus of its normalisation. 
Then $m(C)=l(M_{g,0}(C,[C]))$. If moreover $C$ is rational and 
contained in a 
smooth $K3$ surface $X$, then $e(\jb C)=
l(M_{0,0}(X,d),f)$ (length of the zero-dimensional component 
supported at $f$).}

We now sketch briefly the idea of the proof of Theorem 1.

Let $\cC\to B$ be a semi-universal family of deformations 
of a curve $C$
with planar singularities. We prove that the relative 
compactified Jacobian $\bar 
J \cC$ is smooth; moreover, given any deformation $\cC'\to S$ 
of $C$ with a smooth 
base, $\bar J\cC'$ is smooth if and only if the image of $TS$ 
is transversal in 
$TB$ to the $\delta$-codimensional vector space $V$, the 
support of the tangent 
cone to the $\delta$-constant stratum $B^\delta$.

Assume now $C$ is rational and has $p$ as unique singularity. 
We have to show that 
$e(\bar JC)=m(C,p)$. Choose a one-parameter family $W_t$ 
of smooth 
$\delta$-dimensional subspaces of $B$ such that $0\in W_0$,
 $T_{W_0,0}\cap 
V=\{0\}$, and for general $t$ the intersection 
$W_t\cap B^\delta$ is a set of 
$m(C,p)$ distinct points corresponding to nodal 
curves.

Let $\cC_t\to W_t$ be the induced families. Then 
$\bar J\cC_t$ is a family of 
smooth compact varieties, hence $e(\bar J\cC_t)$ 
does not depend on $t$. Arguing as in 
\cite{Y-Z} and \cite{B}, we prove that 
$e(\bar J\cC_0)=e(\jb C)$, while $e(\jb \cC_t)=m(C,p)$ 
for $t$ general.

{\bf Conventions.} In this paper we will always work 
over the complex numbers, and 
open will mean open in the strong (euclidean) topology 
(unless of course we  specify Zariski open).

{\bf Preliminaries.} We will use the language of 
deformation functors; we recall a 
few facts about them for the reader's convenience.

A deformation functor $D$ will always be a covariant
 functor from local artinian 
$\complex$-algebras to sets, satisfying Schlessinger's 
conditions $(H1)$, $(H2)$, 
$(H3)$, hence admitting a hull (see \cite{Sch}). 
In particular $D$ admits a 
finite-dimensional tangent space, which we denote by 
$TD$, functorial in $D$. A 
functor is smooth if its hull is. The dimension of 
the functor will be equal to 
the dimension of the hull. We will need the following 
elementary result.

{\bf Lemma.} {\sl Let $X\to Y$ and $Z\to Y$ be morphisms 
of smooth deformation 
functors. Then $X\times_YZ$ is smooth of dimension 
$\dim X+\dim Z-\dim Y$ if and 
only if the images of $TX$ and $TZ$ span $TY$.}\\
{\bf Proof.} Base change considerations reduce the 
problem to the case of 
prorepresentable functors, where it is obvious.\qed

It would be possible to replace  deformation functors 
with contravariant  functors 
on the category of germs of complex spaces, and the 
hull with the base of a 
semi-universal family of deformations. 
The two viewpoints correspond to working with formal 
versus convergent power 
series.

{\bf Acknowledgements.} This paper was written at the 
Mittag-Leffler Institute in 
Stockholm, during a special year on Enumerative Geometry. 
The authors are grateful 
for the support received and for making  our 
collaboration possible.
The first author is a member of GNSAGA of CNR.

\section{A. Deformations of curves and sheaves.}

Let $C$ be a reduced projective curve, with 
singular set $\Sigma$.
Any deformation $\cC \lra S$ of $C$ over a 
base $S$ induces a deformation of its
singularities. More
precisely, one can introduce the {\em functor
 of local deformations} by
letting
$\Def^{loc}(C)(T)$ be the set of isomorphism
 classes of data $(U_i,U_i^T)_{i\in 
I}$, where $(U_i)_{i\in I}$ is an affine open 
cover of $C$ and, 
for each $i$, $U_i^T$ is a deformation of $U_i$ 
over $T$; we require that the 
induced deformations of $U_{ij}:=U_i\cap U_j$ 
be the same. 
There is a natural
transformation of functors $\loc: \Def(C) \lra \Def^{loc}(C)$;
 the induced map of 
tangent spaces can be identified with the edge homomorphism
$${\bf T}_C^1 \lra H^0(\cT_C^1)$$
of the local-to-global spectral sequence for the $\cT^i$. 
The kernel of this map is $H^1(\Theta_C)$, the cokernel 
injects in
$H^2(\Theta_C)$ which is zero. The obstruction space 
${\bf T}_C^2$ sits
in an exact sequence
$$0 \lra H^1(\cT_C^1) \lra {\bf T}_C^2 \lra H^0(\cT_C^2) 
\lra 0.$$
As $C$ is reduced, $\cT_C^1$ is supported on a finite 
set of points, hence 
$H^1(\cT_C^1)=0$. If $C$ has locally complete
intersection singularities, then also $\cT_C^2 =0$, so that
in that case ${\bf T}_C^2=0$. Hence in such a situation, and 
in particular when 
$C$ is  a reduced curve with only planar singularities,  the
functors $\Def(C)$ and $\Def^{loc}(C)$  are smooth and $\loc$ 
is a smooth map.

Let $\cF$ be a torsion free coherent sheaf on $C$. Analogously
, we denote by  
$\Def(C,\cF)$ the functor of deformations of the pair, and 
define the functor of 
local deformations by letting $\Def^{loc}(C,\cF)(T)$ be the
 set of isomorphism 
classes of data $(U_i,U_i^T,F_i^T)_{i\in I}$ where
 $(U_i)_{i\in I}$ is an affine 
open cover of $C$, and for each $i$, $(U_i^T,F_i^T)$ 
is a $T$-deformation
of $(U_i,\cF|_{U_i})$ such that the induced deformations 
on $U_{ij}$ are the same.

Again we have a localisation map 
$\Def(C,\cF)\to \Def^{loc}(C,\cF)$.
The four functors introduced sit in a natural commutative diagram
$$
\begin{array}{ccc}
D(C,\cF) & \lra & \Def^{loc}(C,\cF)\\
\downarrow&&\downarrow\\
D(C) & \lra & \Def^{loc}(C)\\
\end{array}
$$
with horizontally localisation maps and vertically forget maps.
Note that this diagram in general is {\em not} cartesian.  
 
{\bf Proposition A.1.} {\sl 
The canonical map
$$ D(C,\cF) \lra D(C) \times_{\Def^{loc}(C)}\Def^{loc}(C,\cF)$$
is smooth.}

{\bf Proof.} 
We have to show the following: Let $\cF_T$, $C_T$ be flat
 deformations of
$C$ and $\cF$ over $T$, $\xi_T \in \Def^{loc}(C,\cF)(T)$ 
the induced local
deformation. If we are given lifts $C_{T'}$ and $\xi_{T'}$
 over a
small extension $T'$ of $T$, then we can lift $\cF_T$ to 
a deformation
$\cF_{T'}$ of $\cF$ over $C_{T'}$ inducing $\xi_{T'}$. 
This can be done
as follows: choose an affine open cover $U_i$ of $C$ 
such that $\xi_{T'}$
is defined by coherent sheaves $F_i'$ on the induced 
cover $U_{i,T'}$ of $C_{T'}$. 
Assume also that $U_{ij}:=U_i\cap U_j$ is smooth for every $i\ne j$.

Let $F_i$ be the restriction of $F_i'$ to $U_{i,T}$. 
The fact that $\cF$ induces 
$\xi_T$ means that we can find isomorphisms 
$\phi_i:\cF_T|_{U_{i,T}}\to F_i$.
The $\phi_i$ induce isomorphisms $\phi_{ij}:F_i\to F_j$ 
over $U_{ij,T}$,  
satisfying the cocycle condition. What we need to prove 
is that the $\phi_{ij}$ 
can be lifted to $\phi_{ij}':F'_i\to F'_j$, again 
satisfying the cocycle  
condition; then the $\phi_{ij}'$ can be used to glue 
together the $F_i'$'s to
a coherent sheaf $\cF_{T'}$ as required. But on $U_{ij}$ 
all the sheaves under 
consideration are line bundles, hence the obstruction to 
the existence of such a 
lifting is an element in $H^2(C_,\cO_{C})$, which is 
zero as $C$ has dimension 
$1$. \qed

If $R$ is a ring and $M$ is an $R$-module, we denote by $\Def(R)$, 
respectively $\Def(R,M)$ the corresponding deformation functors.

{\bf Lemma A.2}. {\sl Let $C$ be a reduced projective curve, $\cF$ 
a torsion free 
module on $C$. Let $\Sigma$ denote the singular locus.
Then the natural morphisms of functors $$
\Def^{loc}(C)\to \prod_{p\in \Sigma}\Def(\cO_{C,p})\qquad
\hbox{\sl and}\qquad
\Def^{loc}(C,\cF)\to \prod_{p\in \Sigma}\Def(\cO_{C,p},\cF_p)$$
are isomorphisms.}

{\bf Proof.} Both morphisms are clearly injective. On the other
 hand, 
surjectivity is obvious since on the smooth open locus, every 
infinitesimal 
deformation is locally trivial and every torsion free sheaf is 
locally free.
\qed

{\bf Proposition A.3}. {\sl Let $P$ be a regular local ring of
 dimension $2$, 
$f\in P$ a nonzero element, and $R=P/(f)$; assume that $R$ is 
reduced. Let $M$ be 
a finitely generated, torsion free
$R$-module of rank $1$. Then $\Def(R,M)$ is a smooth functor.}

{\bf Proof.}
As it is torsion free, the module $M$ has depth $1$. By the 
Auslander-Buchsbaum theorem (see e.g. \cite{Ma}), $M$  has a
free resolution of length $1$ as a $P$-module, so is 
represented as the cokernel 
of some 
$n \times n$ matrix $A$ with entries from $P$:
$$ 0 \lra  P^n \stackrel{A}{\lra}P^n \lra M \lra 0$$

As $M$ is an $R$-module of rank $1$, the determinant ideal 
$(det(A))$ is equal
to $(f)$.\\   
Any flat deformation $M_T$ of $M$ over $T$ (as $P$-module) 
is obtained by 
deforming the matrix 
$A$ to a matrix $A_T$ with entries from 
$P_T:=T \otimes_{\complex} P$, so that
$M_T$ has a presentation
$$0 \lra  P_T^n \stackrel{A_T}{\lra}P_T^n \lra M_T \lra 0.$$
There is a unique deformation $R_T$ of $R$ over $T$
such that $M_T$ is a flat $R_T$-module, given by the ideal 
$(det(A_T))$.
It follows that the natural transformation
$$\Def(A) \lra \Def(R,M)\qquad A_T \mapsto 
(P_T/det(A_T),Coker(A_T))$$   
is {\em smooth}.
As $\Def(A)$, the functor of deformations of the matrix $A$, 
is clearly smooth, 
the functor $\Def(R,M)$ is also smooth. \qed

Note that in the assumption of A.3, although both functors 
$D(R,M)$ and $D(R)$ are 
smooth, the forgetful morphism $D(R,M)\to D(R)$ is not smooth 
in general.

{\bf Remark A.4.} Let $R$ be a one-dimensional local 
${\bf C}$-algebra, and let $M$ be a finitely generated torsion free 
$R$-module. Let $\hat R$ be the completion of $R$, and
 $\hat M=M\otimes_R\hat R$. 
The natural morphism $\Def(R,M)\to \Def(\hat R, \hat M)$
is smooth and induces an isomorphism on tangent spaces, 
and the same is true for 
$\Def(R)\to \Def(\hat R)$. In fact it is easy to see that
 the induced morphisms of 
tangent and obstruction spaces are isomorphisms.

\section{B. Relative compactified Jacobians.}
For any flat projective family of curves $\cC\to S$ 
we let $\jb \cC\to S$ 
be the 
relative compactified Jacobian (see \cite{R}). 
For every closed point $s \in S$ 
the fiber over $s$ of $\jb \cC$ is canonically 
isomorphic to the compactified 
Jacobian $\jb C_s$;
in particular, its points correspond to isomorphism 
classes of torsion free rank 
$1$ degree zero sheaves on $C_s$.

Fix a point $\cF \in \jb \cC$ over $s\in S$, and 
denote again by $(\jb\cC,\cF)$ and $(S,s)$ 
the deformation functors induced by the respective
 germs of complex spaces.
Let $C=C_s$. Remark that if $\cC \lra S$ is a 
semi-universal family of  
deformations of $C$, then we have an isomorphism 
of functors
$$ (\jb\cC, \cF) \simeq  D(C,\cF) $$  

For a general flat family one has a natural 
commutative diagram 
$$
\begin{array}{ccc}
(\jb\cC,\cF) & \lra & \Def^{loc}(C,\cF)\\
\downarrow&&\downarrow\\
(S,s) & \lra & \Def^{loc}(C)\\
\end{array}
$$
and analogous to {\bf A.1.} one has:

{\bf Proposition B.1.} {\sl 
The canonical map
$$ (\jb\cC,\cF) \lra (S,s) \times_{\Def^{loc}(C)}
\Def^{loc}(C,\cF)$$
is smooth.}

We omit the proof, which is almost identical to
 that of {\bf A.1.}

{\bf Corollary B.2.} {\sl 
Let $C$ be a reduced curve with only plane curve
 singularities.
If $\cC\to S$ is a versal family of deformations 
of $C$, then $\jb \cC$ is smooth along $\jb C$, and 
$\jb C$ has local complete intersection singularities.}

{\bf Proof.} The family is versal if and only if the
 natural map $S\to \Def(C)$ is 
smooth. This in turn implies that $S\to \Def^{loc}(C)$ 
is smooth, hence the first 
claim follows from Proposition B.1. On the other hand, 
all fibres of $\jb \cC\to 
S$ have the same dimension $g_a(C)$, therefore each of 
them has local complete 
intersection singularities. \qed

{\bf Corollary B.3.} {\sl With the same assumptions as 
B.2, let $\cC'\to S'$ be 
any deformation of $C$ with smooth base $S'$. Let $\kf$ 
be a torsion free rank $1$ 
degree zero coherent sheaf on $C$. Then the relative 
compactified Jacobian $\jb 
\cC'$ is smooth at $[\kf]$ if and only if the image 
of $TS'$ in $T\Def^{loc}(C)$ 
is transversal to the image of $T\Def^{loc}(C,\kf)$.}

{\bf Proof.} We keep the notation of B.2. The dimension 
of $\bar J\cC'$ is equal 
to $\dim S'+g_a(C)$. 
As $\jb \cC'$ is equal to the fibred product of 
$\jb \cC$ and $S'$ over $S$, 
it follows that $\jb \cC'$ is smooth at $[\kf]$ 
if and only if the image  
of $TS'$ in $TS$ is transversal to that of 
$T(\jb \cC, \kf)$. Proposition B.1 
implies that the image of $T(\jb \cC, \kf)$ 
is the inverse image of the image of
$T\Def^{loc}(C,\kf)$ in $T\Def^{loc}(C)$.\qed

\section{C. The canonical sub-space $V$}

Let $C$ be a reduced curve with only planar 
singularities. 
In this section we study the map 
$$\Def^{loc}(C,\cF)\to \Def^{loc}(C)$$
at the level of tangent spaces. As both functors 
are products corresponding to the 
singularities of $C$ (Lemma A.2) and the tangent 
spaces only depend on the formal 
structure of the singularity (Remark A.4), it suffices 
to analyse what happens for 
$$
\Def(R,M)\to \Def(R)$$
where  $P=\complex[[x,y]]$, $R=P/(f)$, $f$ a non-zero 
element of the maximal ideal 
such that $R$ is reduced, and $M$  a torsion free rank one 
$R$-module given by a presentation 
$$ 0 \lra  P^n \stackrel{A}{\lra}P^n \lra M \lra 0.$$
{\bf Proposition C.1.} {\sl The image of the map $
T\Def(R,M)\to T\Def(R)$
is the image of the first Fitting ideal $F_1(M)$ in the 
quotient ring $T\Def(R)=P/(f,\partial_xf,\partial_yf)$.}\\
{\bf Proof.}
Let $E_{i,j}$ be the $n \times n$ matrix which has entry 
$(i,j)$ equal to $1$ and
all other  entries equal to zero. If $\epsilon^2=0$, then
  $\det(A+\epsilon.E_{i,j})=\det(A)+\epsilon \wedge^{n-1}(A)_{i,j}$,
therefore
we see that by perturbing the matrix $A$ to first order, we 
generate precisely the 
ideal
of $(n-1) \times (n-1)$ minors of the matrix $A$ as  first 
order perturbations of 
$f$.
This is by definition the first Fitting ideal of $F_1(M)$. 
\hfill $\Diamond$

Another description of the ideal $F_1(M)$ is the following

{\bf Proposition C.2.} {\sl $F_1(M)$ is the set of elements
 $r\in R$ such that 
$r=\varphi(m)$ for some $m\in M$, $\varphi\in Hom_R(M,R)$.
}\\
{\bf Proof.} 
As $M$ is maximal Cohen-Macaulay, a resolution of $M$ as an
 $R$-module will be 
$2$-periodic of the form $$
\dots \lra R^n \stackrel{\bar B}{\lra} R^n \stackrel{\bar A}
{\lra} R^n \lra M \lra 
0$$
for some $n\times n$ matrix with $P$-coefficients $B$ with 
the property that $$
AB=BA=f{\bf 1}$$
and $\bar A$, $\bar B$ are the induced matrices with $R$ 
coefficients
(see \cite{E} or \cite{Yo}).

From the $2$-periodicity it follows that there is an exact 
sequence $$
0\lra M \lra R^n \stackrel{\bar A}{\lra} R^n \lra M\lra 0,$$
where $M=\ker A=\hbox{\rm im}\, B$. We split this sequence into $$
\begin{array}{c}
0\lra M\lra R^n \lra N\lra 0\\
0\lra N \lra R^n \lra M \lra 0.\end{array}$$
As $N$ is also torsion free and $R$ is Gorenstein, 
$Ext^1_R(N,R)=0$ by local 
duality. Hence we see from the first sequence that 
the map $Hom_R(R^n,R)\lra 
Hom_R(M,R)$ is surjective.

From this it follows that the ideal obtained by evaluating all 
homomorphisms $\phi\in Hom_R(M,R)$ on all elements of $M$
is the same as the ideal generated by the entries of the 
matrix $\bar B$.

As $M$ has rank $1$, it follows that $det(A)=f$,
and hence the matrix $B$ is the Cramer matrix 
$(\Lambda^{n-1}A)^{tr}$ of $A$. The claim follows. 
\hfill $\Diamond$

Locally, the normalisation $\widetilde C \lra C$ corresponds 
to the inclusion of $R$ in its integral closure $\overline R$
$$R \hookrightarrow \overline R.$$
Recall that the conductor is the ideal $I=Hom_R(\overline R,R)$.
One has 
$$I\subset R\subset \overline R$$
and 
$dim(R/I)=dim(\overline R/R)=\delta(C,p)$.

As an important corollary of {Proposition C.2} we have

{\bf Corollary  C.3.} 
 $F_1(M)\supset I$.\\
{\bf Proof.} Write $\overline R=\oplus \overline R_i$, with
 $\overline R_i$ a 
domain isomorphic to $\complex[[t]]$. Let $Q(\overline R_i)$
 be the quotient field 
of $\overline R_i$, and let $Q(R)=\oplus Q(\overline R_i)$ be 
the total quotient ring of $R$.
As $M$ has rank $1$, $M\otimes_RQ(R)$ is isomorphic to $Q(R)$; 
as it is 
torsion-free, the natural map $M\to M\otimes_RQ(R)$ is 
injective. Hence up to 
isomorphism we can assume that $M$ is a submodule of $Q(R)$. 
Let $m\in M$ be an 
element of minimal valuation
(it exists as $M$ is finitely generated). Then multiplication 
by $m^{-1}$, an 
isomorphism of $Q(R)$ as an $R$-module, sends $M$ to a submodule
 of $\overline  R$ 
containing $1$.

So we can assume that $R\subset M\subset \overline R$.
Let
$c$ be any element of $I$.
Multiplication by $c$ defines a homomorpism
$\phi\in Hom_R(M,R)$ with
$\phi(1)=c$ (note that $1\in R\subset M$). Hence
$$\bigl\{\phi(m)\bigm| m\in M, \ \phi\in Hom(M,R)\bigr\}\supset I.$$
\hfill $\Diamond$

{\bf Remark  C.4.} From the above description one also sees
that $F_1(\overline R)=I$. Hence the differential of the map
 $\Def(R,M)\to 
\Def(R)$ has minimal rank for $M=\overline R$.

Let $C$ be a reduced projective curve with only planar 
singularities, $\Sigma$ its 
singular locus. For $p\in \Sigma$, let $V_p$ be the 
subspace of codimension 
$\delta(C,p)$ in 
$TD(C,p)$ generated by the conductor, and put
$$V^{loc}=\prod_{p \in \Sigma} V_p \subset 
T\Def^{loc}(C) =\prod_{p \in 
\Sigma}TD(C,p).$$ 
Let $V$ be the inverse image of $V^{loc}$ in $T\Def(C)$; note 
that $V$ is a linear 
subspace of codimension $\delta(C)$.
If $B$ is the base space of a semi-universal family of deformations 
of $C$, then $TB$ is identified with $TD(C)$.

{\bf Proposition C.5.} {\sl Let $\cC\to B$ be a semi-universal 
family of  
deformations of $C$. Then for any $\cF\in \bar JC$ the image of
 the tangent map 
$\jb \cC\to B$ at $\cF$ contains the subspace $V$, and there
 exists at least
one such $\cF$ for which the image is exactly $V$. }
\\
{\bf Proof.} The first statement follows immediately from 
Proposition C.1 and 
Corollary C.3, by applying Proposition B.1 and Lemma A.2. 
The second statement 
follows in the same way from Remark C.4; e.g., we can take
 $\cF=n_*(\cO_{\tilde 
C})$, where $n:\tilde C\to C$ is the normalisation map.
\qed

\section{D. The $\delta$-constant stratum.}

Let $C$ be a reduced curve with only planar singularities. 
We denote by $B$  an 
appropriate representative of the semi-universal
deformation of $C$. The stratum $B^\delta$ is defined as 
the set of points where 
the geometric genus of the fibres is constant.
This amounts to saying that 
$$\sum_{x\in C_t} \delta(C_t,x)$$
is constant for $t\in B^\delta$
and equal to $g_a(C)-g(\widetilde C)$, hence the name.

The analytic set $B^\delta$ (we give it the reduced induced 
structure) is very 
singular in general, but its properties 
can be related directly to the local $\delta$-constant strata 
$$B^\delta(C,p).$$
To be more precise, $B^\delta$ is the pullback of
$B^{\delta,loc}=\prod B^\delta(C,p)$ under the smooth map
$B\lra B^{loc}$.
So let $(C,p)\subset (\complex^2,0)$ be a reduced  plane curve 
singularity, with normalisation
$$(\widetilde C,q)\stackrel{n}{\lra}(C,p),\qquad q=n^{-1}(p).$$
Note that in general $q$ will be a finite set of distinct points,
 one for each 
branch of $C$ at $p$.
We denote for brevity by $\Def(n)$ the functor of deformations of 
$n:(\widetilde 
C,q)\to (C,p)$ (that is, we are allowed to deform $C$ and 
$\tilde C$ as well as 
the map).

{\bf Lemma D.1.} {\sl $D(n)$ is smooth.}\\
{\bf Proof.} The  morphism $D((C,p)\to (\complex^2,0))
\lra D(n)$ (given by taking the image of the deformation of the map) 
is smooth.
Hence it is enough to verify that $D((C,p)\to (\complex^2,0)$ is smooth,
and this is obvious.\qed

{{\bf Theorem D.2.} (\cite{T}, \cite{D-H}).}
{\sl Let $(C,p)\subset (\complex^2,0)$  be a reduced plane curve 
singularity,  
$n:(\widetilde C,q) \lra(C,p) $ its normalisation.
Let $B(C,p)$ be a semi-universal family for $\Def(C,p)$ and
$$B^{\delta}(C,p)\subset B(C,p)$$
the $\delta$-constant stratum. Then one has:\\
{(1)} 
The normalisation $\widetilde B^{\delta}(C,p)$ of 
$B^\delta(C,p)$ is a smooth space. \\
{(2)} The pullback of the semi-universal family to 
$\tilde B^\delta$ admits a 
simultaneous resolution of singularities. This makes
$\widetilde B^{\delta}(C,p)$ into a semi-universal
family for $\Def(n)$.\\
{(3)}
The codimension of $B^\delta\subset B$ is $\delta(C,p)$.
Over the generic point $p\in B^{\delta}$, the curve $C_p$ 
has precisely 
$\delta(C,p)$ double points as its only singularities.\\
(4) The tangent cone to the $\delta$-constant stratum is 
supported on $V_p$, the 
vector subspace generated by the conductor ideal. \qed}

The second half of (2) is in fact not explicitly stated 
in either of \cite{T}, 
\cite{D-H}; however it follows easily from Lemma D.1. A 
similar argument is 
presented in the proof of Proposition F.2, so we don't 
repeat it here.

\section{E.  Proof of Theorem 1.}

Let $C$ be a reduced projective rational curve with only planar 
singularities.
We want to show that $e(\bar JC)=m(C)$. In particular let  
$(C,p)$ be a reduced 
plane curve singularity. Let $C$
be a projective rational curve that has $(C,p)$ as its only 
singular point.
Then it follows that  $e(\bar JC)=m(C,p)$.

Let $\Phi:\cC\lra B$ be a semi-universal family of  
deformations of $C$;
we denote its fibres by $C_s=\Phi^{-1}(s)$, $C_0=C$.
Let $\pi:\jb\cC \lra B$ be the corresponding family of compactified
Jacobians.
We always assume that we have chosen discs as representatives 
for the corresponding germs.
We may also assume that the induced morphism
$j:B\lra B^{loc}$ is smooth and has contractible fibres.
We choose a section $\sigma:B^{loc}\lra B$ of $j$ with $\sigma(0)=0$.
We will denote $\overline B:=\sigma(B^{loc})$, 
$\overline B^\delta:=\sigma(B^\delta)$ and $\overline V:=\sigma(V)$.

Let $(W,0)\subset (\overline B,0)$ be a smooth subspace of dimension
$\delta+1$ containing the point $(0,0)$ together with a smooth map
$\lambda:(W,0)\lra (T,0)$ to a disc $(T,0)\subset (\C,0)$.
$W$ is a one parameter family of $\delta$-dimensional subspaces 
$W_t=\lambda^{-1}(t)\subset \overline B$.
We require in addition that $W_0$ is transverse to $V$.

\bigskip
        \epsfxsize 4cm
        \epsfysize 4cm
        \centerline{\epsfbox{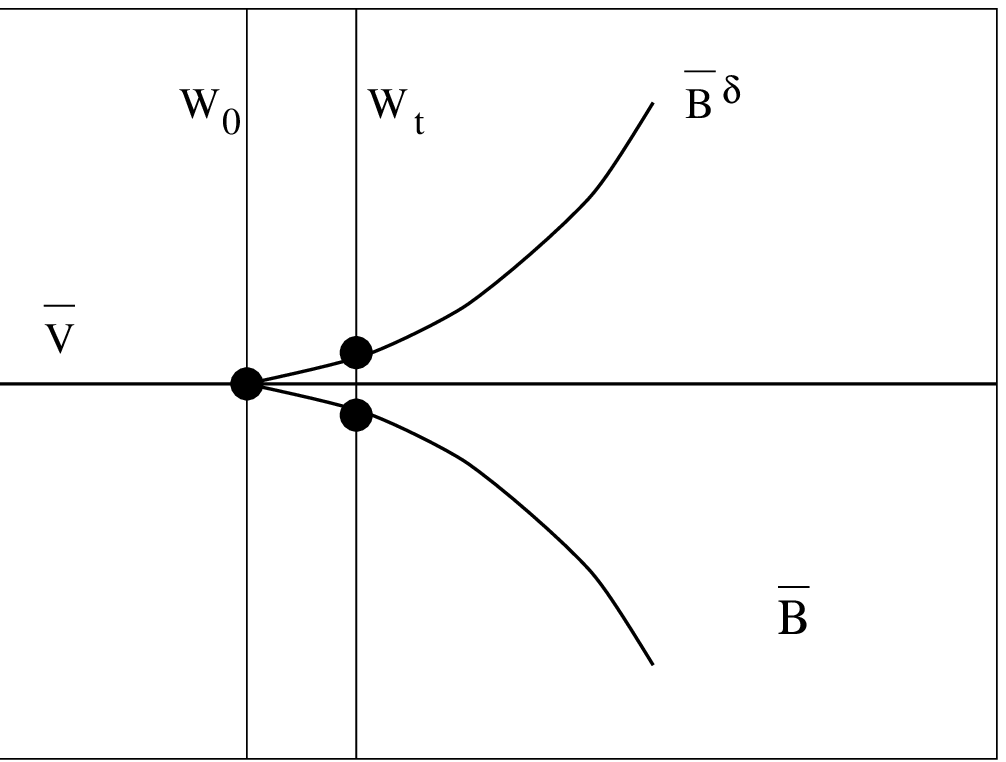}}
\bigskip

 By {Theorem D.2} we can choose $W$ in such a way that 
for $t\ne 0$ the fibre $W_t$ intersects $\overline B^\delta$ in 
$mult(B^{\delta})$ points, and for $s\in W_t\cap \overline B^\delta$ the  
corresponding curve $C_s$ has precisely $\delta$ nodes as singularities.
For $s\in  W_t \setminus \overline B^\delta$ the curve $C_s$ will have
positive genus. Let $\bar \Delta\subset B$ be a closed disc, and let $Z
=W\cap \bar\Delta$.
We define the family $\rho:\jb\cC_Z\lra T$ by pullback:

\def\mapd#1{\Big\downarrow 
\rlap{$\vcenter{\hbox{$\scriptstyle#1$}}$}}

 $$\begin{array}{ccccc}
&&\jb\cC_Z&\lra&\jb\cC\\
&\stackrel{\rho}{\swarrow}&\mapd{\pi}&&\mapd{\pi}\\
T&\stackrel{\lambda}{\longleftarrow}&Z&\lra&B
\end{array}
$$
As we have chosen $W_0$ to be transversal to $V$, Proposition C.5
implies that $\rho$ is smooth along $\pi^{-1}(0)$; by making 
$\bar\Delta$ and $T$ 
smaller we can assume
that $\rho$ is smooth. As $\rho$ is also proper,  
all the fibres $\rho^{-1}(t)$ are diffeomorphic, in particular 
they all have the 
same Euler number.

The space $\rho^{-1}(t)$ is the union, for $s\in W_t$, of
 $\bar JC_s$. We know 
that if $C_s$ has positive geometric genus, then $e(\bar JC_s)$ 
is zero; arguing 
as in \cite{B}, we obtain that $$
e(\rho^{-1}(t))=\sum_{s\in W_t\cap \bar B^\delta}e(\bar JC_s)$$
(note that if $s\in W_t$, then $C_s$ is rational if and only 
if $s\in \bar 
B^\delta$).\\
The intersection of $W_0\subset \overline B$ with 
$\overline B^\delta$
consists only of the point $0$ corresponding to the curve $C$.
Therefore $e(\rho^{-1}(0))=e(\jb C)$.\\
On the other hand, for $t\ne 0$, $W_t$ intersects 
$\overline B^{\delta}$ in $mult( B^\delta)$ points and for $s\in 
\overline B^{\delta}\cap W_t$ the curve $C_s$ has precisely 
$\delta$
nodes as singularities.
As for a nodal rational curve $C_s$, the Euler number 
$e(\jb C_s)$ is equal to $1$, we obtain
$$e(\rho^{-1}(t))=
\sum_{s\in W_t}e(\jb C_s)=\sum_{s\in W_t\cap \overline B^{\delta}}1=
mult(B^\delta).$$
So we get
$$e(\jb C)=e(\rho^{-1}(0))=e(\rho^{-1}(t))=mult(B^\delta).$$
\hfill$\Diamond$

\section{F. The invariant as length of moduli of stable maps}

Let $C$ be a reduced projective curve with only plane curve 
singularities; let 
$n:\tilde C\to C$ be its normalisation, and $g$ the genus of 
$\tilde C$.  
Let $m(C)=\prod m(C,p)$. The scheme $\overline M_{g,0}(C,[C])$ 
parametrizing 
stable birational maps from a genus $g$ curve to $C$ contains 
only one point, 
namely the normalisation of $C$.
The aim of this section is to prove that its length is equal 
to $m(C)$.
Note that if $C$ is an isolated rational curve inside a 
smooth manifold $Y$, 
$\overline M_{g,0}(C,[C])$ is naturally a closed subscheme
 of $\overline  
M_{g,0}(Y,[C])$; in particular, $m(C)$ is a lower bound 
for the length of the 
corresponding component of $M_{g,0}(Y,[C])$ (in case this 
scheme also has 
dimension zero).

Denote by $\Def(n)$ the deformation functor of the triple 
$(n:\widetilde C\to C)$, 
and by $\Def^{loc}(n)$ the corresponding local deformation 
functor. As before 
$\Def^{loc}(n)$ is the product over the singular points $p$ 
of $C$ of $\Def(n,p)$, 
the deformation functor of the triple $n:(\tilde C, 
n^{-1}(p))\to (C,p)$. 

If $(C,p)$ is the germ of a planar reduced curve singularity, 
then $\Def(n,p)$ is 
a smooth functor (see section D). 

{\bf Lemma F.1.} {\sl
The natural morphism of functors $\Def(n)\to 
\Def^{loc}(n)\times_{\Def^{loc}(C)}\Def(C)$ is an isomorphism.}

{\bf Proof.} Let $C_T$ be an infinitesimal deformation of $C$, 
and let $U_i$ be an 
open cover of $C$ such that $U_{ij}$ is smooth for each $i\ne j$. 
Let $V_i=n^{-1}(U_i)$. Let $U_{i,T}$ be the deformation of $U_i$ 
induced by $C_T$, 
and assume we are given a deformation $n_{i,T}:V_{i,T}\to U_{i,T}$
 of 
$n_i:=n|_{V_i}$. Then to lift $(C_T,n_{i,T})$ to a deformation 
of $n$ we must 
choose gluing isomorphisms $\psi_{ij}:V_{ij,T}\to V_{ji,T}$ 
satisfying the cocycle 
condition and compatible with the other data, namely the maps 
$n_{i,T}$ and the 
gluing isomorphisms $\phi_{ij}:U_{ij,T}\to U_{ji,T}$ induced 
by $C_T$. But 
$U_{ij}$ is smooth, so that $n|_{V_{ij}}$ is an isomorphism 
for each $i\ne j$; 
hence the $\psi_{ij}$ are univocally determined by the 
$\phi_{ij}$ and 
automatically satisfy the cocycle condition.
\qed

Let us now denote by $B(\cdot)$ the germ of complex
 space being a hull for the 
functor $\Def(\cdot)$. Note that Lemma F.1 implies 
that there is a cartesian 
diagram $$
\cdia{B(n)}{B(C)}{B^{loc}(n)}{B^{loc}(C).}$$

{\bf Proposition F.2.} {\sl Let $C$ be a reduced projective  
curve with planar 
singularities, $n:\tilde C\to C$ be the normalisation, 
$g=g(\tilde C)$.
Let $\pi:{\cal C}\to B(C)$ be a semi-universal deformation
 of $C$. 
Denote by $M= M_{g,0}({\cal C},[C])$; then $M$ is smooth at 
$n$, and the natural 
map $M\to B^\delta:=B^\delta(C)$ is the normalisation map.}

{\bf Proof.} Write $M$ for the germ of $M$ at $n$.
As the domain of $n$ is a smooth curve, the same is true 
for all stable maps in a 
neighborhood of $n$. Hence $M$ is  isomorphic to $B(n)$. 
By Lemma F.1, together 
with Lemma D.1,   we deduce that $B(n)$ is smooth. By the 
definition of $B^\delta$  
the natural map $M\to B(C)$ factors via $B^\delta$, hence, 
as $M$ is smooth, via 
its normalisation $\tilde B^\delta$. On the other hand, we
 know that the family 
$\tilde {\cal C}\to \tilde B^\delta$ gotten by pullback 
admits a very weak 
simultaneous resolution of singularities \cite{T}, inducing a morphism 
$\tilde B^\delta\to M$. It is easy to check pointwise that 
these two morphisms are 
inverse to each other (both $\tilde B^\delta$ and $M$
just parametrize the normalisation maps of the fibres of $\pi$). 
As 
both $\tilde B^\delta$ and $M$ are smooth, a bijective morphism 
must be an 
isomorphism. \qed

{\bf Proof of Theorem 2}. The scheme $M_{g,0}(C,[C])$ is 
the fibre over the point 
$[C]$ of the morphism $\tilde B^\delta\to B^\delta$; this 
is the multiplicity of 
$B^\delta$ at $[C]$, as $\tilde B^\delta$ is smooth. This 
proves the first 
equality.

Let now $X$ be a smooth projective surface, $C\subset X$ a 
reduced irreducible 
curve, $n:\tilde C\to C$ the normalisation, $g=g(\tilde C)$.
 Assume that 
$n$ is an isolated point of $\overline{M}_{g,0}(X,[C])$, and 
let $M_n$ be the 
connected component of $n$. $M_n$ contains $M_{g,0}(C,[C])$ 
as a closed subscheme, 
so we always have an inequality $$
l(M_n)\ge l(M_{g,0}(C,[C]))=m(C).$$
This inequality is an equality if and only the natural morphism 
$M_n\to Hilb(X)$ 
sending each map to its image factors scheme theoretically
 (and not only  
set-theoretically) via $C$.

Hence to complete the proof of Theorem 2, it is enough to 
show that this is the 
case if $C$ is rational and $X$ is a $K3$ surface. Let $S$ 
be the complete linear 
system defined by $C$ on $X$, and let $\cC\to S$ be the 
universal curve.
It is known that  $\bar J\cC$ is smooth, see \cite{Mu}; but 
this means precisely that $S$ maps transverse 
to the $\delta$-constant stratum in $B(C)$, and we are done 
in view of Corollary 
B.3.\qed

\section{G. Examples.}

\def\el{\ell}
\def\barM{M}
\def\alp{\alpha}
\def\AAA{{\bf A}}
\def\P{{\bf P}}

{\bf Example 1 (Beauville):} {\sl Let $(C,o)$ be the 
singularity of equation 
$x^q=y^p$, with $p<q$ and $(p,q)=1$. Then $$m(C,o)=
{{1}\over{p+q}}{p+q\choose p}.$$}

{\bf Proof.} We write for simplicity $\overline M(X,\beta)$ 
instead of $\overline 
M_{0,0}(X,\beta)$; if $X$ is a curve and $\beta=[X]$ we omit 
it. Let $C$ be the 
plane curve of equation $y^pz^{q-p}=x^q$. $C$ is a rational 
curve with two 
singular points, $o=(0,0,1)$ and $\infty=(1,0,0)$.
Let $\alp:C'\to C$ be the partial normalisation of $C$ at 
$\infty$. 
By Theorem 2, it is enough to prove that 
$$l(\overline M(C'))={{1}\over{p+q}}{p+q\choose p}=:N(p,q).$$
The natural map $\overline M(C')\to \overline M(C)$ given by 
$\mu\mapsto \alp\circ 
\mu$ is a closed embedding, and the closed subscheme 
$\overline M(C')$ is 
identified by requiring the deformation of the normalisation 
morphism to be 
locally trivial near $\infty$.
On the other hand $\overline M(C)$ is naturally a closed subset
 of $\overline 
M(\P^2,q\el)$, where $\el$ is the class of a line. 

Let $n:\P^1\to C$ be the normalisation map, and choose coordinates 
on $\P^1$ such 
that $n(s,t)=(t^ps^{q-p},t^q,s^q)$. A morphism in $\overline 
M(\P^2,q\el)$ near 
$n$ has equations $$(t^ps^{q-p}+x,t^q+y,s^q+z),$$ for suitable
 homogeneous 
polynomials $x,y,z$ of degree $q$. 

We  impose  the conditions that the image of the map be 
contained in $C$
and  that the deformation be locally trivial at $\infty$. 
Then we   eliminate the indeterminacy generated by
a reparametrization of $\P^1$ and a rescaling of the 
coordinates on $\P^2$. 
We get that all deformations of $n$ in $\overline M(C)$ 
must be (in affine 
coordinates where $z=1$) of the form $$
t\mapsto (t^p+\textstyle\sum\limits_{i=0}^p x_i t^i,t^q+
\textstyle\sum\limits_{i=0}^qy_it^i).$$

Hence we are now left with the following problem: 
compute the length
of the $\complex$--algebra with generators 
$x_0,\ldots,x_{p-2},y_0,\ldots,y_{q-2}$ 
and relations given by the coefficients of the polynomial 
$f^q-g^p$, where 
$f=t^p+\sum x_it^i$ and $g=t^q+\sum y_it^i$.

It is easy to check that the equation $f^q=g^p$ is 
equivalent to $qf'g=pg'f$
by taking $d/dt\circ\log$ on both sides. The $t$-degree 
of $qf'g-pg'f$ is $p+q-1$, 
however we only get $p+q-2$ equations as the coefficients 
of $t^{p+q-1}$ and 
$t^{p+q-2}$ are zero anyway.
Moreover, if we consider the variables $x_i$ (resp.\ $y_i$) 
as having degree 
$p-i$ (resp.\ $q-i$), the equations we obtain are homogeneous of
degree $2,\ldots,p+q-1$.

Now we recall the weighted B\'ezout theorem, which says that
 if we have a
zero-dimensional algebra given by $N$ homogeneous equations of 
degrees $e_j$ in $N$ weighted variables of degree $d_j$, then
 the length of the 
algebra is $\prod e_j/\prod d_j$.

Applying the formula in our case, with $N=p+q-2$, 
$(d_j)=(2,3,\ldots,p,2,3,\ldots,q)$ and 
$e_j=(2,3,\ldots,p+q-1)$ gives
$$
N(p,q)={\prod e_j\over \prod d_j}={(p+q-1)!\over p!q!}
={1\over p+q}{p+q\choose 
p}.$$

{\bf Example 2.} We would like to outline an algorithm for
 the computation of
$m(C,p)$ for a planar, reduced and irreducible curve 
singularity $(C,p)$. 
Assume we know how to realize $(C,p)$ as singularity 
of a rational curve. It is 
then easy to realize it as singularity of a plane 
rational curve $C$, whose other 
singularities are only nodes. Let $d$ be the degree 
of the curve, $F(x,y,z)=0$ its 
equation, and $\bar n=(\bar x,\bar y,\bar z)$ an explicit
 normalisation given by  
homogeneous polynomials of degree $d$ in $s,t$. Assume 
without loss of generality 
that $\bar z$ contains the monomial $s^d$ with nonzero 
coefficient.

Then we can describe the scheme $M_{0,0}(C,[C])$ explicitly 
 as follows.  
Choose three points $p_i$ ($i=1,2,3$) in $\IP^1$ mapping via 
$n$ to smooth points 
of $C$; let $L_i\subset \IP^2$ be a line transversal to $C$ 
at $n(p_i)$.

Choose variables $x_i$, $y_i$ and $z_i$ for $i=0,\ldots d$, 
and
let $x$ be the polynomial $\bar x+\sum_i x_i s^it^{d-i}$; 
define $y$ and $z$ in a 
similar way.

Then $M_{0,0}(C,[C])$ is naturally isomorphic the 
subscheme of $
\hbox{\sl Spec}\, \complex [x_i,y_i,z_i]$ defined 
by the equations 
$$\begin{array}{c}
z_d=0,\\
(x,y,z)(p_i)\in L_i \qquad i=1,2,3,\\
F(x,y,z)=0.
\end{array}$$
In fact, all deformations of $\bar n$ are again morphisms
 of degree $d$ from 
$\IP^1$ to $\IP^2$, hence are given by polynomials of 
degree $d$.  
The first four equations, defining a linear subspace,
correspond to choosing local coordinates near $\bar n$ 
on $M_{0,0}(\P^2,d)$;   the 
last one, which is a system of 
$d^2$ equations, imposes the condition that the 
scheme-theoretic image of the 
morphism be contained in 
$C$.\\

Mittag-Leffler Institute, Aurav\"agen 17, S 18262 Djursholm, Sweden.\\
E-mail: {\tt fantechi\@@science.unitn.it, lothar\@@mpim-bonn.mpg.de,}\\
 {\tt straten\@@mathematik.uni-mainz.de}
\end{document}